\input amstex
\documentstyle{amsppt}
\topmatter \magnification=\magstep1 \pagewidth{5.2 in}
\pageheight{6.7 in}
\abovedisplayskip=10pt \belowdisplayskip=10pt
\parskip=8pt
\parindent=5mm
\baselineskip=2pt
\title
 A note on multiple Dirichlet's $q$-$L$-function
\endtitle
\author  Taekyun Kim   \endauthor

\affil{ {\it Institute of Science Education,\\
        Kongju National University, Kongju 314-701, S. Korea}\\
          {\it e-mail: {\rm tkim$\@$kongju.ac.kr }}\\\\}\endaffil
\thanks  2000 Mathematics Subject Classification:  11S80, 11B68, 11M99 .\endthanks
\abstract{ Recently, the two variable $q$-$L$-functions which
interpolate the generalized $q$-Bernoulli polynomials associated
with $\chi$ are introduced and studied, cf. [2]. In this paper, we
construct multiple Dirichlet's $q$-$L$-function which interpolates
the generalized multiple $q$-Bernoulli numbers attached to $\chi $
at negative integer and study its properties.  }\endabstract

\leftheadtext{Taekyun Kim}
\endtopmatter

\document

\head \S 1. Introduction \endhead
 Let $0<q<1$ be and  for any positive integer $k$, define its
 $q$-analogue $[k]_q=\frac{1-q^k}{1-q}$. Let $\Bbb C$ be the field
 of complex numbers.
The $q$-Bernoulli numbers are usually defined  as
$$\beta_{0,q}=\frac{q-1}{\log q}, \text{ } \left(q\beta_q+1
\right)^n-\beta_{n,q}=\delta_{n,1}, \tag1$$ where $\delta_{n,1}$
is Kronecker symbol and we use the usual convention about
replacing $\beta_q^i$ by $\beta_{i,q}, $ cf. [2, 4, 5, 6]. Note
that $\lim_{q\rightarrow 1}\beta_{k,q}=B_k$, where $B_k$ are the
$k$-th ordinary Bernoulli numbers. In [3, 5, 6], the $q$-Bernoulli
polynomials are also defined by
$$F_q(t,x)=\frac{q-1}{\log
q}e^{\frac{t}{1-q}}-t\sum_{n=0}^{\infty}q^{n+x}e^{[n+x]_qt}
=\sum_{n=0}^{\infty}\frac{\beta_{n,q}(x)}{n!}t^n, \text{ for
$x\in\Bbb C$}. \tag 2$$ From (1) and (2), we can derive the below
formula:
$$\beta_{n,q}(x)=\sum_{i=0}^n\binom ni
q^{xi}\beta_{i,q}[x]_q^{n-i}=\left(\frac{1}{1-q}\right)^n\sum_{i=0}^n\binom
ni(-1)^i q^{xi}\frac{i}{[i]_q},$$ where $\binom ni$ is the
binomial coefficient.
 In the recent paper the multiple
generalized Bernoulli numbers attached to $\chi$,
$B_{n,\chi}^{(r)},$ are defined by
$$F_{\chi}^r(t) =\sum_{a_1,\cdots,a_r=1}^f\frac{\chi(\sum_{i=1}^r
a_i)t^re^{(\sum_{i=1}^ra_i)t}}{(e^{ft}-1)^r}=\sum_{n=0}^{\infty}B_{n,\chi}^{(r)}\frac{t^n}{n!},
\text{ for $|t|<2\pi/f$, (see [1]). } \tag3 $$ In this paper we
consider the $q$-analogue of the above multiple generalized
Bernoulli numbers attached to $\chi$. Finally we construct the
multiple $q-L$-function which interpolates the $q$-analogue of
multiple generalized Bernoulli numbers attached to $\chi$ at
negative integers and investigate its properties.

 \head 2. Multiple Dirichlet's $q-L$-function    \endhead

Let $\chi (\neq 1)$ be the primitive Dirichlet's  character with
conductor $f\in\Bbb N$. Then the generalized $q$-Bernoulli numbers
attached to $\chi$ are defined as
$$F_{q,\chi}(t)=-t\sum_{a=1}^f\chi(a)\sum_{n=0}^{\infty}q^{fn+a}e^{[fn+a]_qt}=\sum_{n=0}^{\infty}
\beta_{n,\chi,q}\frac{t^n}{n!}, \text{ cf. [5].} \tag 4$$ By
(2)and (4), we easily see that
$$\beta_{n,\chi,q}=[f]_q^{n-1}\sum_{a=1}^f
\chi(a)\beta_{n,q^f}(\frac{a}{f}).$$ By the meaning of Eq.(4), we
can consider the multiple generalized $q$-Bernoulli numbers
attached to $\chi$ as follows:
$$\aligned
F_{q,\chi}^r (t)&=(-t)^r\sum_{a_1,\cdots,
a_r=1}^f\chi(\sum_{i=1}^r a_i)\sum_{n_1,\cdots,
n_r=0}^{\infty}q^{\sum_{i=1}^r (a_i+n_if)}
e^{([\sum_{i=1}^r(a_i+n_if)]_q)t}\\
&=\sum_{n=0}^{\infty}\beta_{n,\chi,q}^{(r)}\frac{t^n}{n!}.\endaligned
\tag5$$ Note that $\lim_{q\rightarrow
1}\beta_{n,\chi,q}^{(r)}=B_{n,\chi}^{(r)}.$ Let $\Gamma (s)$ be
the gamma function. Then it is easy to see that
$$\aligned
&\frac{1}{\Gamma (s)}\int_{0}^{\infty}F_{\chi,q}^r(-t)t^{s-1-r}dt\\
&=\sum_{n_1.\cdots,n_r=1}^{\infty}\chi(n_1+\cdots+n_r)q^{n_1+\cdots+n_r}\frac{1}{\Gamma
(s)}\int_{0}^{\infty}t^{s-1}e^{-[n_1+\cdots+n_r]_qt}dt\\
&=\sum_{n_1,\cdots,
n_r=1}^{\infty}\frac{\chi(n_1+\cdots+n_r)}{[n_1+\cdots+n_r]_q^s}q^{n_1+\cdots+n_r}.
\endaligned \tag6$$
Thus, we can define the multiple Dirichlet's $q$-$L$-function as
follows:
$$L_q^r(s,\chi)=\sum_{n_1,\cdots,n_r=1}^{\infty}q^{n_1+\cdots+n_r}\frac{\chi(n_1+\cdots+n_r)}{[n_1+\cdots+n_r]_q^s},
\text{ for $s\in\Bbb C$, }\tag7$$ where $\chi$ is non-trivial
primitive Dirichlet's character with conductor $f\in\Bbb N .$ By
(5), (6) and (7), we easily see that
$$L_q^r(-n,\chi)=(-1)^r\frac{n!}{(n+r)!}\beta_{n+r,\chi,q}^{(r)},
\text{ for $n\in\Bbb N$}. \tag8$$ Let $s$ be a complex variable,
$a$ and $F$ be integers with $0<a<F$. We now consider the function
$H_{r,q}(s;a_1,\cdots,a_r|F)$ as follows:
$$H_{r,q}(s;a_1,\cdots, a_r|F)=\sum_{\Sb m_1,\cdots,m_r>0\\
m_i\equiv a_i (\mod
F)\endSb}\frac{q^{m_1+\cdots+m_r}}{[m_1+\cdots+m_r]_q^s}=[F]_q^{-s}\zeta_{r,q^F}(s,\frac{a_1+\cdots+a_r}{F}),\tag9$$
where
$\zeta_{r,q}(s,a)=\sum_{n_1,\cdots,n_r=0}^{\infty}\frac{q^{n_1+\cdots+n_r+a}}{[n_1+\cdots+n_r+a]_q^s}
, \text{ (see [5, 6]). } $

The function $H_{r,q}(s;a_1,\cdots, a_r|F)$ is a meromorphic for
$s\in\Bbb C $ with poles at $s=1,\cdots, r.$ In [5], the multiple
$q$-Bernoulli polynomials are defined by
$$(-t)^r\sum_{n_1,\cdots,n_r=0}^{\infty}q^{x+\sum_{i=1}^r
n_i}e^{[x+\sum_{i=1}^rn_i]_qt}=\sum_{n=0}^{\infty}B_{n,q}^{(r)}(x)\frac{t^n}{n!}.
\tag10 $$ For $n\in\Bbb N ,$ it was well known that
$$\zeta_{r,q}(-n,x)=(-1)^r\frac{n!}{(n+r)!}B_{n+r,q}^{(r)}(x),
\text{ cf. [5]. }\tag11 $$ Let $F\in\Bbb N$ be the conductor of
$\chi .$ Then the multiple Drichlet's $q$-$L$-function can be
expressed as the sum
$$L_q^r(s,\chi)=\sum_{a_1,\cdots,a_r=1}^F\chi(a_1+\cdots+a_r)H_{r,q}(s;a_1,\cdots,a_r|F),
\text{ for $s\in\Bbb C$.} \tag 12$$ By (9) and (12), we easily see
that
$$H_{r,q}(-n;a_1,\cdots,a_r|F)=[F]_q^n(-1)^r\frac{n!}{(n+r)!}B_{n+r,
q^F}^{(r)}(\frac{a_1+\cdots+a_r}{F}), \text{ for $n\in\Bbb
N.$}\tag13$$ From (5) and (10), we can also derive the below
formula:
$$\beta_{n,\chi,q}^{(r)}=[F]_q^n\sum_{a_1,\cdots,a_r=1}^F
\chi(a_1+\cdots+a_r)B_{n,q^F}^{(r)}(\frac{a_1+\cdots+a_r}{F}).\tag14$$
By using (12), (13) and (14), we easily see that
$$L_q^r(-n,\chi)=(-1)^r\frac{n!}{(n+r)!}\beta_{n+r,\chi,q}^{(r)}.$$
In Eq.(10), we note that
$$B_{n,q}^{(r)}(x)=\sum_{k=0}^n\binom nk
[x]_q^{n-k}q^{xk}B_{k,q}^{(r)}, \text { cf. [5, 6]}, \tag15$$
where $B_{n,q}^{(r)}=B_{n,q}^{(r)}(0).$ From Eq.(15), we derive
the function $H_{r,q}(s;a_1,\cdots,a_r|F)$ which is modified by
$$\aligned
&H_{r,q}(s;a_1,\cdots,a_r|F) \\
&=\frac{1}{[F]_q^r}\frac{[\sum_{i=1}^r
a_i]_q^{-s+r}}{\prod_{j=1}^r (s-j)}
\sum_{k=0}^{\infty}\binom{-s+r}{k}\left(\frac{[F]_q}{[\sum_{i=1}^ra_i]_q}\right)^kq^{(\sum_{i=1}^ra_i)k}B_{k,q}^{(r)}.
\endaligned\tag16$$
By (12) and (16), we obtain the following:
$$\aligned
L_q^r(s,\chi)&=\frac{1}{\prod_{j=1}^r
(s-j)}\frac{1}{[F]_q^r}\sum_{a_1,\cdots,a_r=1}^F\chi(a_1+\cdots+a_r)[a_1+\cdots+a_r]_q^{-s+r}\\
&\cdot\sum_{m=0}^{\infty}\binom{r-s}{m}\left(\frac{[F]_q}{[a_1+\cdots+a_r]_q}
\right)^m q^{(a_1+\cdots+a_r)m}B_{m,q}^{(r)}.
\endaligned\tag17$$
Finally, we suggest the below question

\proclaim{ Question}  Is it possible to give the $p$-adic analogue
of  Eq.(17) which can be viewed as interpolating, in the same way
that $L_{p,q}(s, \chi)$ interpolates $L_q(s, \chi) \text{  in [2,
7 ] } ?$
\endproclaim

\enddocument